 \documentclass[11pt]{article}
 \usepackage{mathrsfs}

 \usepackage{mathrsfs,amsfonts}
 \usepackage[dvips]{color}

 \newtheorem{theorem}{Theorem}[section]
 \newtheorem{lemma}[theorem]{Lemma}
 \newtheorem{corol}[theorem]{Corollary}
 \newtheorem{prop}[theorem]{Proposition}
 \newtheorem{example}[theorem]{Example}

 \def\blemma{\begin{lemma}\sl{}\def\elemma{\end{lemma}}}
 \def\bproposition{\begin{prop}\sl{}\def\eproposition{\end{prop}}}
 \def\btheorem{\begin{theorem}\sl{}\def\etheorem{\end{theorem}}}
 \def\bcorollary{\begin{corol}\sl{}\def\ecorollary{\end{corol}}}
 \def\bexample{\begin{example}\rm{}\def\eexample{\end{example}}}

 \def\bitemize{\begin{itemize}}\def\eitemize{\end{itemize}}
 \def\itm{\item}

 \def\beqlb{\begin{eqnarray}}\def\eeqlb{\end{eqnarray}}
 \def\beqnn{\begin{eqnarray*}}\def\eeqnn{\end{eqnarray*}}

 \def\proof{\noindent{\it Proof.~~}}\def\qed{\hfill$\Box$\medskip}

 \def\<{\langle}\def\>{\rangle}

 \def\mcr{\mathscr}\def\mbb{\mathbb}\def\mbf{\mathbf}

 \def\ar{\!\!&}\def\nnm{\nonumber}

 \def\itDelta{{\it\Delta}}

 \def\lline{---------------------------------}

\begin{document}

\noindent{\ }

\bigskip\bigskip\bigskip

\centerline{\LARGE\bf Strong solutions of jump-type}

\medskip

\centerline{\LARGE\bf stochastic equations\footnote{ Supported by NSFC
(No.~11131003), 973 Program (No.~2011CB808001) and 985 Program.}}

\bigskip\bigskip

\centerline{Zenghu Li and Fei Pu}

\bigskip

\centerline{School of Mathematical Sciences, Beijing Normal University,}

\centerline{Beijing 100875, People's Republic of China}

\centerline{E-mails: {\tt lizh@bnu.edu.cn} and {\tt
pufei@mail.bnu.edu.cn}}

\bigskip

\centerline{\lline\lline\lline}

\medskip

{\narrower

\noindent\textit{Abstract.} We establish the existence and uniqueness of
strong solutions to some jump-type stochastic equations under
non-Lipschitz conditions. The results improve those of Fu and Li
\cite{FuL10} and Li and Mytnik \cite{LiM11}.

\smallskip


\smallskip


\smallskip

\noindent\textbf{Abbreviated Title:} Strong solutions of stochastic
equations

\par}

\bigskip\bigskip

\centerline{\lline\lline\lline}

\bigskip


\section{Introduction}

The problem of existence and uniqueness of solutions to jump-type
stochastic equations under non-Lipschitz conditions have been studied by
many authors; see, e.g., \cite{Bas03, Bas04, BBC04, Fou11, FuL10, Kom82,
LiM11} and the references therein. In particular, some criteria for the
existence and pathwise uniqueness of non-negative and general solutions
were given in \cite{Fou11, FuL10, LiM11}. Stochastic equations have played
important roles in the recent progresses in the study of continuous-state
branching processes; see, e.g., \cite{BeL06, DaL06, DaL12, Li12}. The main
difficulty of pathwise uniqueness for jump-type stochastic equations
usually comes from the compensated Poisson integral term. Let us consider
the equation
 \beqlb\label{1.1}
dx(t) = \phi(x(t-))d\tilde{N}(t),
 \eeqlb
where $\{\tilde{N}(t): t\ge 0\}$ is a compensated Poisson process. For
each $0<\alpha<1$ there is a $\alpha$-H\"older continuous function $\phi$
so that the pathwise uniqueness for (\ref{1.1}) fails. In fact, before the
first jump of the Poisson process, the above equation reduces to
 \beqlb\label{1.2}
dx(t) = -\phi(x(t))dt.
 \eeqlb
Then to assure the pathwise uniqueness for (\ref{1.1}) the uniqueness of
solution for (\ref{1.2}) is necessary. If we set $h_\alpha(x) =
(1-\alpha)^{-1}x^\alpha 1_{\{x\ge 0\}}$, then both $x_1(t)=0$ and $x_2(t)
= t^{1/(1-\alpha)}$ are solutions of (\ref{1.2}) with $\phi = -h_\alpha$.
From those it is easy to construct two distinct solutions of (\ref{1.1}).
The key of the pathwise uniqueness results in \cite{FuL10, LiM11} is to
consider a non-decreasing kernel for the compensated Poisson integral term
in the stochastic equation. The condition was weakened considerably by
Fournier \cite{Fou11} for stable driving noses. In fact, as a consequence
of Theorem~2.2 in \cite{LiM11}, given any $x(0)\in \mbb{R}$ there is a
pathwise unique strong solution to (\ref{1.1}) with $\phi = h_\alpha$. On
the other hand, the monotonicity assumption also excludes some interesting
jump-type stochastic equations. Two of them are given below.

\bexample\label{t1.1} Let $z^2\nu(dz)$ be a finite measure on $(0,1]$.
Suppose that $\tilde{M}(ds,dz,dr)$ is a compensated Poisson random measure
on $(0,\infty)\times (0,1]^2$ with intensity $ds\nu(dz)dr$. Given $0\le
x(0)\le 1$, we consider the stochastic integral equation
 \beqlb\label{1.3}
x(t) = x(0) + \int_0^t\int_0^1\int_0^1 zq(x(s-),r) \tilde{M}(ds,dz,dr),
 \eeqlb
where
 \beqnn
q(x,r) = 1_{\{r\le 1\land x\}} - (1\land x)1_{\{x\ge 0\}}.
 \eeqnn
This equation was introduced by Bertoin and Le~Gall \cite{BeL05} in their
study of generalized Fleming-Viot flows. The existence and uniqueness of a
weak solution flow to (\ref{1.3}) was proved in \cite{BeL05}. The pathwise
uniqueness for the equation follows from a result in \cite{DaL12}. The
result cannot be derived directly from the those in \cite{FuL10, LiM11}
since $x\mapsto q(x,r)$ is not a non-decreasing function. \eexample

\bexample\label{t1.2} Let $(1\wedge u^2) \mu(du)$ be a finite measure on
$(0,\infty)$. Suppose that $\tilde{N}(ds,du,dr)$ is a compensated Poisson
random measure on $(0,\infty)^3$ with intensity $ds\mu(du)dr$. Given
$y(0)\ge 0$, we consider the stochastic equation
 \beqlb\label{1.4}
y(t) = y(0) + \int_0^t\int_0^\infty\int_0^\infty g(y(s-),u,r)
\tilde{N}(ds,du,dr),
 \eeqlb
where
 \beqnn
g(x,u,r) = - 1_{\{rx\le 1\}}x(1-e^{-u}).
 \eeqnn
Some generalizations of the above equation were introduced by D\"{o}ring
and Barczy \cite{DoB11} in the study of self-similar Markov processes.
From their results it follows that (\ref{1.4}) has a pathwise unique
non-negative strong solution. Since $x\mapsto g(x,u,r)$ is not
non-decreasing, one cannot derive the pathwise uniqueness for (\ref{1.4})
from the results in \cite{FuL10, LiM11}. \eexample

In this paper, we give some criteria for the existence and pathwise
uniqueness of strong solutions of jump-type stochastic equations. The
results improve those in \cite{FuL10, LiM11} and can be applied to
equations like (\ref{1.3}) and (\ref{1.4}). In Section~2 we give some
basic formulations of the stochastic equations. Two theorems on the
pathwise uniqueness of general solutions are presented in Section~3. In
Section~4 we prove the existence of weak solutions by a martingale problem
approach. The main results on the existence and pathwise uniqueness of
general strong solutions are given in Section~5. In Section~6 we give some
results on the existence and pathwise uniqueness of non-negative strong
solutions. Throughout this paper, we make the conventions
 \beqnn
\int_a^b = \int_{(a,b]}
 \quad\mbox{and}\quad
\int_a^\infty = \int_{(a,\infty)}
 \eeqnn
for any $b\ge a\ge 0$. Given a function $f$ defined on a subset of
$\mbb{R}$, we write
 \beqnn
\itDelta_zf(x) = f(x+z) - f(x)
 \quad\mbox{and}\quad
D_zf(x) = \itDelta_zf(x) - f^\prime(x)z
 \eeqnn
if the right hand sides are meaningful.


\section{Preliminaries}

\setcounter{equation}{0}

Suppose that $\mu_0(du)$ and $\mu_1(du)$ are $\sigma$-finite measures on
the complete separable metric spaces $U_0$ and $U_1$, respectively.
Throughout this paper, we consider a set of parameters
$(\sigma,b,g_0,g_1)$ satisfying the following basic properties:
 \bitemize

\itm $x\mapsto \sigma(x)$ is a continuous function on $\mbb{R}$;

\itm $x\mapsto b(x)$ is a continuous function on $\mbb{R}$ having the
decomposition $b=b_1-b_2$ with $b_2$ being continuous and non-decreasing;

\itm $(x,u)\mapsto g_0(x,u)$ and $(x,u)\mapsto g_1(x,u)$ are Borel
functions on $\mbb{R}\times U_0$ and $\mbb{R}\times U_1$, respectively.

 \eitemize

Let $(\Omega, \mcr{G}, \mcr{G}_t, \mbf{P})$ be a filtered probability
space satisfying the usual hypotheses. Let $\{B(t): t\ge 0\}$ be a
standard $(\mcr{G}_t)$-Brownian motion and let $\{p_0(t): t\ge 0\}$ and
$\{p_1(t): t\ge 0\}$ be $(\mcr{G}_t)$-Poisson point processes on $U_0$ and
$U_1$ with characteristic measures $\mu_0(du)$ and $\mu_1(du)$,
respectively. Suppose that $\{B(t)\}$, $\{p_0(t)\}$ and $\{p_1(t)\}$ are
independent of each other. Let $N_0(ds,du)$ and $N_1(ds,du)$ be the
Poisson random measures associated with $\{p_0(t)\}$ and $\{p_1(t)\}$,
respectively. Let $\tilde{N}_0(ds,du)$ be the compensated measure of
$N_0(ds,du)$. By a \textit{solution to} the stochastic equation
 \beqlb\label{2.1}
x(t)
 \ar=\ar
x(0) + \int_0^t \sigma(x(s-))dB(s) + \int_0^t\int_{U_0} g_0(x(s-),u)
\tilde{N}_0(ds,du) \nnm \\
 \ar\ar
+ \int_0^t b(x(s-))ds + \int_0^t\int_{U_1} g_1(x(s-),u) N_1(ds,du),
 \eeqlb
we mean a c\`adl\`ag and $({\mcr{G}}_t)$-adapted real process $\{x(t)\}$
that satisfies the equation almost surely for every $t\ge0$. Since $x(s-)
\neq x(s)$ for at most countably many $s\ge0$, we can also use $x(s)$
instead of $x(s-)$ for the integrals with respect to $dB(s)$ and $ds$ on
the right hand side of (\ref{2.1}). We say \textit{pathwise uniqueness}
holds for (\ref{2.1}) if for any two solutions $\{x_1(t)\}$ and
$\{x_2(t)\}$ of the equation satisfying $x_1(0) = x_2(0)$ we have $x_1(t)
= x_2(t)$ almost surely for every $t\ge 0$. Let $(\mcr{F}_t)_{t\ge 0}$ be
the augmented natural filtration generated by $\{B(t)\}$, $\{p_0(t)\}$ and
$\{p_1(t)\}$. A solution $\{x(t)\}$ of (\ref{2.1}) is called a
\textit{strong solution} if it is adapted with respect to $(\mcr{F}_t)$;
see~\cite[p.163]{IkW89} or \cite[p.76]{Sit05}. Let $U_2\subset U_1$ be a
set satisfying $\mu_1(U_1\setminus U_2)<\infty$. We also consider the
equation
 \beqlb\label{2.2}
x(t)
 \ar=\ar
x(0) + \int_0^t \sigma(x(s-))dB(s) + \int_0^t\int_{U_0} g_0(x(s-),u)
\tilde{N}_0(ds,du) \nnm \\
 \ar\ar
+ \int_0^t b(x(s-))ds + \int_0^t\int_{U_2} g_1(x(s-),u) N_1(ds,du).
 \eeqlb

\bproposition\label{t2.1} If (\ref{2.2}) has a strong solution for every
given $x(0)$, so does (\ref{2.1}). If the pathwise uniqueness holds for
(\ref{2.2}), it also holds for (\ref{2.1}). \eproposition

The above proposition can be proved similarly as Proposition~2.2
in~\cite{FuL10}. Then all conditions in the paper only involve $U_2$
instead of $U_1$.


\section{Pathwise uniqueness}

\setcounter{equation}{0}

In this section, we prove some results on the pathwise uniqueness for
(\ref{2.2}) under non-Lipschitz conditions. Suppose that
$(\sigma,b,g_0,g_1)$ are given as in the second section. Let us consider
the following conditions on the modulus of continuity:
 \bitemize

\itm[{\rm(3.a)}] for each integer $m\ge 1$ there is a non-decreasing and
concave function $z\mapsto r_m(z)$ on $\mbb{R}_+$ such that $\int_{0+}
r_m(z)^{-1}\, dz=\infty$ and
 \begin{eqnarray*}
|b_1(x)-b_1(y)| + \int_{U_2} |l_1(x,y,u)| \mu_1(du)\le r_m(|x-y|),
\qquad |x|,|y|\le m,
 \end{eqnarray*}
where $l_1(x,y,u) = g_1(x,u)-g_1(y,u)$;

\itm[{\rm(3.b)}] the function $x\mapsto x + g_0(x,u)$ is non-decreasing
for all $u\in U_0$ and for each integer $m\ge 1$ there is a constant
$K_m\ge 0$ such that
 \beqnn
|\sigma(x)-\sigma(y)|^2 + \int_{U_0} l_0(x,y,u)^2 \mu_0(du)\le K_m|x-y|,
\qquad |x|,|y|\le m,
 \eeqnn
where $l_0(x,y,u) = g_0(x,u) - g_0(y,u)$.

 \eitemize

Let us define a sequence of functions $\{\phi_k\}$ as follows. For each
integer $k\ge 0$ define $a_k = \exp\{-k(k+1)/2\}$. Then $a_k\to 0$
decreasingly as $k\to \infty$ and
 \beqnn
\int_{a_k}^{a_{k-1}}z^{-1} dz = k, \qquad k\ge 1.
 \eeqnn
Let $x\mapsto \psi_k(x)$ be a non-negative continuous function supported
by $(a_k,a_{k-1})$ so that
 \beqlb\label{3.1}
\int_{a_k}^{a_{k-1}}\psi_k(x)dx=1
 \quad\mbox{and}\quad
\psi_k(x)\le 2(kx)^{-1}
 \eeqlb
for every $a_k< x< a_{k-1}$. For $z\in \mbb{R}$ let
 \beqlb\label{3.2}
\phi_k(z)=\int_0^{|z|}dy\int_0^y\psi_k(x)dx.
 \eeqlb
It is easy to see that the sequence $\{\phi_k\}$ has the following
properties:
 \bitemize

\itm[{\rm(i)}] $\phi_k(z)\mapsto |z|$ non-decreasingly as $k\rightarrow
\infty$;

\itm[{\rm(ii)}] $0\leq \phi_k^\prime(z)\leq 1$ for $z\geq 0$ and $-1\leq
\phi_k^\prime(z)\leq 0$ for $z\leq 0$;

\itm[{\rm(iii)}] $0\le |z|\phi_k^{\prime\prime}(z) = |z|\psi_k(|z|)\le
2k^{-1}$ for $z\in \mbb{R}$.

 \eitemize
By Taylor's expansion, for any $h,\zeta\in \mbb{R}$ we have
 \beqlb\label{3.3}
D_h\phi_k(\zeta)
 =
h^2\int_0^1\psi_k(|\zeta+th|)(1-t)dt
 \le
\frac{2}{k}h^2\int_0^1 {(1-t)\over |\zeta+th|}dt.
 \eeqlb

\blemma\label{t3.1} Suppose that $x\mapsto x + g_0(x,u)$ is non-decreasing
for $u\in U_0$. Then, for any $x\neq y\in \mbb{R}$,
 \beqlb\label{3.4}
D_{l_0(x,y,u)}\phi_k(x-y)
 \le
\frac{2}{k}\int_0^1 {l_0(x,y,u)^2(1-t)\over |x-y + tl_0(x,y,u)|}dt
 \le
{2l_0(x,y,u)^2\over k|x-y|}.
 \eeqlb
\elemma

\proof The first inequality follows from (\ref{3.3}). Since $x\mapsto x +
g_0(x,u)$ is non-decreasing, for $x>y\in \mbb{R}$ we have
$x-y+l_0(x,y,u)\ge 0$, and hence $x-y + tl_0(x,y,u)\ge 0$ for $0\le t\le
1$. It is elementary to see
 \beqnn
\ar\ar\int_0^1 {l_0(x,y,u)^2(1-t)\over x-y + tl_0(x,y,u)}dt \cr
 \ar\ar\qquad
= l_0(x,y,u)\int_0^1 \Big[{x-y + l_0(x,y,u)\over x-y + tl_0(x,y,u)} -
1\Big]dt \cr
 \ar\ar\qquad
= \big[x-y + l_0(x,y,u)\big] \log\Big(1 + {l_0(x,y,u)\over x-y}\Big) -
l_0(x,y,u) \cr
 \ar\ar\qquad
\le [x-y + l_0(x,y,u)]{l_0(x,y,u)\over x-y} - l_0(x,y,u) \cr
 \ar\ar\qquad
= {l_0(x,y,u)^2\over x-y}.
 \eeqnn
Then the second inequality in (\ref{3.4}) follows by symmetry. \qed

\btheorem\label{t3.2} Suppose that conditions (3.a,b) are satisfied. Then
the pathwise uniqueness for (\ref{2.2}) holds. \etheorem

\proof By condition (3.b) and Lemma~\ref{t3.1}, for $x\neq y\in \mbb{R}$
satisfying $|x|,|y|\le m$ we have
 \beqnn
\phi_k^{\prime\prime}(x-y)[\sigma(x)-\sigma(y)]^2
 \le
K_m\phi_k^{\prime\prime}(x-y)|x-y|
 \le
\frac{2K_m}{k}
 \eeqnn
and
 \beqnn
\int_{U_0} D_{l_0(x,y,u)}\phi_k(x-y) \mu_0(du)
 \le
\int_{U_0} {2l_0(x,y,u)^2\over k|x-y|} \mu_0(du)
 \le
\frac{2K_m}{k}.
 \eeqnn
The right-hand sides of both inequalities tend to zero uniformly on
$|x|,|y|\le m$ as $k\to \infty$. Then the pathwise uniqueness for
(\ref{2.2}) follows by a simple modification of Proposition~3.1 in
\cite{LiM11}; see also Theorem~3.1 in \cite{FuL10}. \qed

We next introduce some condition that is particularly useful in
applications to stochastic equations driven by L\'evy processes. The
condition is given as follows:
 \bitemize

\itm[{\rm(3.c)}] there is a constant $0\le c\le 1$ such that $x\mapsto cx
+ g_0(x,u)$ is non-decreasing for all $u\in U_0$ and for each integer
$m\ge 1$ there are constants $K_m\ge 0$ and $p_m> 0$ such that
 \beqnn
|\sigma(x)-\sigma(y)|^2\le K_m|x-y| \quad\mbox{and}\quad |l_0(x,y,u)|\le
|x-y|^{p_m}f_m(u)
 \eeqnn
for $|x|, |y|\le m$, where $l_0(x,y,u) = g_0(x,u)-g_0(y,u)$ and $u\mapsto
f_m(u)$ is a strictly positive function on $U_0$ satisfying
 \beqnn
\int_{U_0} [f_m(u)\land f_m(u)^2] \mu_0(du)<\infty.
 \eeqnn

 \eitemize
For each $m\ge 1$ and the function $f_m$ specified in (3.c) we define the
constant
 \beqlb\label{3.5}
\alpha_m = \inf\Big\{\beta>1: ~ \lim_{x\to 0+}  x^{\beta-1}\! \int_{U_0}
f_m(u)1_{\{f_m(u)\ge x\}}\mu_0(du) = 0\Big\}.
 \eeqlb
By Lemma~2.1 in \cite{LiM11} we have $1\le \alpha_m\le 2$.

\blemma\label{t3.3} Suppose that condition (3.c) holds. Then for any $h\ge
0$ and $|x|, |y|\le m$ we have
 \beqnn
\ar\ar\int_{U_0}D_{l_0(x,y,u)}\phi_k(x-y)\mu_0(du) \nnm\\
 \ar\ar\qquad
\le \frac{2}{k}|x-y|^{2p_m-1}1_{\{(1-c)|x-y|<
a_{k-1}\}}\int_{U_0}f_m(u)^21_{\{f_m(u)\le h\}}\mu_0(du)\nnm\\
 \ar\ar\qquad\quad
+\, 2|x-y|^{p_m}1_{\{(1-c)|x-y|< a_{k-1}\}}\int_{U_0}f_m(u)1_{\{f_m(u)>
h\}}\mu_0(du).
 \eeqnn
\elemma

\proof We first consider $x>y\in \mbb{R}$. Since $x \mapsto cx + g_0(x,u)$
is non-decreasing, we have $c(x-y) + l_0(x,y,u)\ge 0$, and hence $c(x-y) +
tl_0(x,y,u)\ge 0$ for $0\le t\le 1$. It follows that $x-y + tl_0(x,y,u)\ge
(1-c)(x-y)$ for $0\le t\le 1$. Then $(1-c)(x-y)\ge a_{k-1}$ implies $x-y +
tl_0(x,y,u)\ge a_{k-1}$ for $0\le t\le 1$. In view of the equality in
(\ref{3.3}) we have
 \beqnn
D_{l_0(x,y,u)}\phi_k(x-y) = 0 \quad\mbox{if}\quad (1-c)(x-y)\ge a_{k-1}.
 \eeqnn
By the symmetry of $\phi_k$ it is follows that, for arbitrary $x,y\in
\mbb{R}$,
 \beqlb\label{3.6}
D_{l_0(x,y,u)}\phi_k(x-y) = 0 \quad\mbox{if}\quad (1-c)|x-y|\ge a_{k-1}.
 \eeqlb
Then we can use condition (3.c) to get
 \beqnn
D_{l_0(x,y,u)}\phi_k(x-y)
 \ar\le\ar
2|l_0(x,y,u)|1_{\{(1-c)(x-y)< a_{k-1}\}} \nnm\\
 \ar\le\ar
2|x-y|^{p_m}f_m(u)1_{\{(1-c)|x-y|< a_{k-1}\}}.
 \eeqnn
Similarly, by (\ref{3.4}) we have
 \beqnn
D_{l_0(x,y,u)}\phi_k(x-y)
 \ar\le\ar
\frac{2l_0(x,y,u)^2}{k|x-y|}1_{\{(1-c)|x-y|<
a_{k-1}\}} \nnm\\
 \ar\le\ar
\frac{2}{k}|x-y|^{2p_m-1}f_m(u)^21_{\{(1-c)|x-y|< a_{k-1}\}}.
 \eeqnn
Those give the desired result. \qed

\btheorem\label{t3.4} Suppose that conditions (3.a,c) hold with: (i)~$c=1,
\alpha_m=2, p_m=1/2$; or (ii)~$c<1, \alpha_m<2, 1-1/\alpha_m< p_m\le 1/2$.
Then the pathwise uniqueness holds for (\ref{2.2}). \etheorem

\proof Let us consider the case (i). By Lemma~\ref{t3.3}, for any $h\ge 1$
and $|x|, |y|\le m$ we have
 \beqnn
\ar\ar\int_{U_0} D_{l_0(x,y,u)}\phi_k(x-y) \mu_0(du) \cr
 \ar\ar\qquad
\le \frac{2}{k}\int_{U_0} f_m(u)^21_{\{f_m(u)\le h\}} \mu_0(du)
+ 2\sqrt{2m}\int_{U_0} f_m(u)1_{\{f_m(u)>h\}} \mu_0(du)  \\
 \ar\ar\qquad
\le \frac{2h}{k}\int_{U_0} [f_m(u)\land f_m(u)^2] \mu_0(du) +
2\sqrt{2m}\int_{U_0} f_m(u)1_{\{f_m(u)>h\}} \mu_0(du).
 \eeqnn
By letting $k\to \infty$ and $h\to \infty$ one can see
 \beqnn
\lim_{k\to \infty}\int_{U_0} D_{l_0(x,y,u)}\phi_k(x-y) \mu_0(du) = 0.
 \eeqnn
Then the pathwise uniqueness for (\ref{2.2}) follows by a modification of
Proposition~3.1 in \cite{LiM11}; see also Theorem~3.1 in \cite{FuL10}. The
case (ii) follows as in the proof of Proposition~3.3 in \cite{LiM11}. \qed

We remark that our conditions (3.b) and (3.c) improve similar conditions
in \cite{FuL10, LiM11}, where it was assumed that $x\mapsto g_0(x,u)$ is
non-decreasing for all $u\in U_0$. The following example shows that the
global monotonicity of the functions $x\mapsto x + g_0(x,u)$ and $x\mapsto
cx + g_0(x,u)$ in conditions (3.b) and (3.c) are necessary to assure the
pathwise uniqueness.

\bexample\label{t3.5} Let us consider the equation (\ref{1.1}). Let $0<
\alpha< 1$ be a constant and define the bounded positive $\alpha$-H\"older
continuous function
 \beqlb\label{3.7}
\phi(x) = (1-\alpha)^{-1}(|x|^\alpha\land |x-1|^\alpha) 1_{\{0\le x\le
1\}}, \qquad x\in \mbb{R}.
 \eeqlb
Clearly, this function is nondecreasing in the interval $(-\infty,1/2)$
and nonincreasing in the interval $(1/2,\infty)$. Let $y_1(t) = 1$ for
$t\ge0$ and let
 \beqnn
y_2(t)
 =
\left\{\begin{array}{ll}
1-t^{1/(1-\alpha)} \ar\mbox{ for $0\le t< 2^{\alpha-1}$,} \\
(2^\alpha-t)^{1/(1-\alpha)} \ar\mbox{ for $2^{\alpha-1}\le t< 2^\alpha$,} \\
0 \ar\mbox{ for $t\ge 2^\alpha$.}
\end{array}\right.
 \eeqnn
It is elementary to show that both $\{y_1(t)\}$ and $\{y_2(t)\}$ are
solutions of (\ref{1.2}) satisfying $y_1(0) = y_2(0) = 1$. Based on
$\{y_1(t)\}$ and $\{y_2(t)\}$, it is easy to construct infinitely many
solutions of (\ref{1.2}) satisfying $y(0)=1$. Therefore (\ref{1.1}) has
infinitely many solutions $\{x(t)\}$ satisfying $x(0)=1$.  \eexample


\section{Weak solutions}

\setcounter{equation}{0}

In this section, we prove the existence of the weak solution to
(\ref{2.2}) by considering the corresponding martingale problem. Let
$(\sigma,b,g_0,g_1)$ be given as in the second section. Let $C^2(\mbb{R})$
be the set of twice continuously differentiable functions on $\mbb{R}$
which together with their derivatives up to the second order are bounded.
For $x\in \mbb{R}$ and $f\in C^2(\mbb{R})$ we define
 \beqlb\label{4.1}
Af(x) \ar=\ar \frac{1}{2}\sigma(x)^2 f^{\prime\prime}(x) +
\int_{U_0}D_{g_0(x,u)}f(x) \mu_0(du) \cr
 \ar\ar
+\, b(x)f^\prime(x) + \int_{U_2} \itDelta_{g_1(x,u)}f(x) \mu_1(du).
 \eeqlb
To simplify the statements we introduce the following condition:
 \bitemize

\itm[{\rm(4.a)}] there is a constant $K\ge 0$ such that
 \beqnn
\ar\ar |b(x)| + \sigma(x)^2 + \int_{U_0} g_0(x,u)^2 \mu_0(du) \cr
 \ar\ar\qquad
+ \int_{U_2} \big[|g_1(x,u)|\vee g_1(x,u)^2\big] \mu_1(du) \le K, \qquad
x\in \mbb{R}.
 \eeqnn

 \eitemize

\bproposition\label{t4.1} Suppose that condition (4.a) holds. Then a
c\`adl\`ag process $\{x(t): t\ge 0\}$ is a weak solution to (\ref{2.2}) if
and only if for every $f\in C^2(\mbb{R})$,
 \beqlb\label{4.2}
f(x(t)) - f(x(0)) - \int_0^tAf(x(s))ds, \qquad t\ge 0
 \eeqlb
is a locally bounded martingale. \eproposition

\proof Without loss of generality, we assume $x(0)\in \mbb{R}$ is
deterministic. If $\{x(t): t\ge 0\}$ is a solution to (\ref{2.2}), by
It\^o's formula it is easy to see that (\ref{4.2}) is a locally bounded
martingale. Conversely, suppose that (\ref{4.2}) is a martingale for every
$f\in C^2(\mbb{R}_+)$. By a standard stopping time argument, we have
 \beqnn
x(t) = x(0) + \int_0^tb(x(s-))ds + \int_0^tds\int_{U_2} g_1(x(s-),u)
\mu_1(du) + M(t)
 \eeqnn
for a square-integrable martingale $\{M(t): t\ge 0\}$. As in the proof of
Proposition~4.2 in \cite{FuL10}, we obtain the equation (\ref{2.2}) on an
extension of the probability space by applying martingale representation
theorems; see, e.g., \cite[p.90 and p.93]{IkW89}. \qed

Now suppose that conditions (3.a,b) and (4.a) are satisfied. For
simplicity, in the sequel we assume the initial value $x(0)\in \mbb{R}$ is
deterministic. Let $\{V_n\}$ be a non-decreasing sequence of Borel subsets
of $U_0$ so that $\cup_{n = 1}^{\infty}V_n = U_0$ and $\mu_0(V_n) <
\infty$ for every $n\ge 1$. It is easy to see that
 \beqnn
x\mapsto \int_{V_n} g_0(x,u) \mu_0(du)
 \eeqnn
is a bounded continuous function on $\mbb{R}$. For $n\ge 1$ and $x\in
\mbb{R}$ let
 \beqlb\label{4.3}
\chi_n(x)=\Bigg\{\begin{array}{ll} n, &\mbox{if}~ x > n, \cr x, \quad
&\mbox{if}~ |x| \leq n, \cr -n, &\mbox{if}~ x < -n.
\end{array}
 \eeqlb
By the result on continuous-type stochastic equations, there is a weak
solution to
 \beqlb\label{4.4}
x(t) \ar=\ar x(0) + \int_0^t\sigma(x(s))dB(s) + \int_0^t b(x(s)) ds \cr
 \ar\ar
- \int_0^tds\int_{V_n} g_0(\chi_n(x(s)),u) \mu_0(du);
 \eeqlb
see, e.g., \cite[p.169]{IkW89}. We can rewrite (\ref{4.4}) into
 \beqlb\label{4.5}
x(t) \ar=\ar x(0) + \int_0^t \sigma(x(s))dB(s) + \int_0^t [b_1(x(s)) +
\mu_0(V_n)\chi_n(x(s))] ds \cr
 \ar\ar
- \int_0^t \Big\{b_2(x(s)) + \int_{V_n} [\chi_n(x(s)) +
g_0(\chi_n(x(s)),u)]\mu_0(du)\Big\} ds,
 \eeqlb
where
 \beqnn
x \mapsto b_2(x) + \int_{V_n} [\chi_n(x) + g_0(\chi_n(x),u)] \mu_0(du)
 \eeqnn
is a bounded continuous non-decreasing function on $\mbb{R}$. By
Theorem~\ref{t3.2} the pathwise uniqueness holds for (\ref{4.5}), so it
also holds for (\ref{4.4}). Then there is a pathwise unique strong
solution to (\ref{4.4}). Let $\{W_n\}$ be a non-decreasing sequence of
Borel subsets of $U_2$ so that $\cup_{n=1}^\infty W_n = U_2$ and
$\mu_1(W_n)<\infty$ for every $n\ge 1$. Following the proof of
Proposition~2.2 in \cite{FuL10} one can see for every integer $n\ge 1$
there is a strong solution to
 \beqlb\label{4.6}
x(t)
 \ar=\ar
x(0) + \int_0^t \sigma(x(s))dB(s) + \int_0^t b(x(s))ds  \cr
 \ar\ar
+ \int_0^t\int_{V_n} g_0(\chi_n(x(s-)),u) \tilde{N}_0(ds,du)  \cr
 \ar\ar
+ \int_0^t\int_{W_n} g_1(x(s-),u) N_1(ds,du).
 \eeqlb
By Theorem~\ref{t3.2} the pathwise uniqueness holds for (\ref{4.6}), so
the equation has a unique strong solution; see, e.g., \cite[p.104]{Sit05}.
Let us denote the strong solution to (\ref{4.6}) by $\{x_n(t): t\ge 0\}$.
By Proposition~\ref{t4.1}, for every $f\in C^2(\mbb{R})$,
 \beqlb\label{4.7}
f(x_n(t)) = f(x_n(0)) + \int_0^tA_nf(x_n(s))ds + \mbox{mart.},
 \eeqlb
where
 \beqnn
A_nf(x) \ar=\ar \frac{1}{2}\sigma(x)^2f^{\prime\prime}(x) + \int_{V_n}
D_{g_0(\chi_n(x),u)}f(x) \mu_0(du)  \cr
 \ar\ar
+\, b(x)f^\prime(x) + \int_{W_n} \itDelta_{g_1(x,u)}f(x) \mu_1(du).
 \eeqnn

\blemma\label{t4.2} Suppose that conditions (4.a) and (3.a,b) are
satisfied. If $x_n\to x$ as $n\to \infty$, then $A_nf(x_n)\to Af(x)$ as
$n\to \infty$. \elemma

\proof Let $M\ge 0$ be a constant so that $|x|, |x_n|\le M$ for all $n\ge
1$. Under the conditions, it is easy to see that
 \beqnn
x\mapsto \int_{V_k^c} g_0(x,u)^2 \mu_0(du) + \int_{W_k^c} |g_1(x,u)|
\mu_1(du)
 \eeqnn
is a continuous function for each $k\ge 1$. By Dini's theorem we have, as
$k\to \infty$,
 \beqnn
\varepsilon_k := \sup_{|x|\le M}\Big[\int_{V_k^c} g_0(x,u)^2 \mu_0(du) +
\int_{W_k^c} |g_1(x,u)| \mu_1(du)\Big]\to 0.
 \eeqnn
Let $y_n=\chi_n(x_n)$. For $n\ge k$ we have
 \beqlb\label{4.8}
\ar\ar\Big|\int_{V_n} D_{g_0(y_n,u)}f(x_n) \mu_0(du) - \int_{U_0}
D_{g_0(x,u)}f(x) \mu_0(du)\Big| \cr
 \ar\ar\qquad
\le \int_{V_k} \Big|D_{g_0(y_n,u)}f(x_n) - D_{g_0(x,u)}f(x)\Big| \mu_0(du)
+ \|f^{\prime\prime}\|\varepsilon_k \cr
 \ar\ar\qquad
\le \int_{V_k} \Big|f(x_n+g_0(y_n,u)) - f(x+g_0(x,u))\Big| \mu_0(du) \cr
 \ar\ar\qquad\quad
+ \int_{V_k} |f(x_n)-f(x)| \mu_0(du) + \|f^{\prime\prime}\|\varepsilon_k
\cr
 \ar\ar\qquad\quad
+ \int_{V_k} \Big|f^\prime(x_n)g_0(y_n,u) - f^\prime(x)g_0(x,u)\Big|
\mu_0(du) \cr
 \ar\ar\qquad
\le \|f^\prime\|\int_{V_k} \Big|(x_n+g_0(y_n,u))-(x+g_0(x,u))\Big|
\mu_0(du) \cr
 \ar\ar\qquad\quad
+ \int_{V_k} |f(x_n)-f(x)| \mu_0(du) + \|f^{\prime\prime}\|\varepsilon_k
\cr
 \ar\ar\qquad\quad
+ \|f^\prime\|\int_{V_k} |g_0(y_n,u)-g_0(x,u)| \mu_0(du) \cr
 \ar\ar\qquad\quad
+ \int_{V_k} |f^\prime(x_n)-f^\prime(x)||g_0(x,u)| \mu_0(du) \cr
 \ar\ar\qquad
\le 2\|f^\prime\|\int_{V_k} |g_0(y_n,u)-g_0(x,u)| \mu_0(du) \cr
 \ar\ar\qquad\quad
+ \Big[\|f^\prime\||x_n-x| + |f(x_n)-f(x)|\Big]\mu_0(V_k) +
\|f^{\prime\prime}\|\varepsilon_k \cr
 \ar\ar\qquad\quad
+ |f^\prime(x_n)-f^\prime(x)|\mu_0(V_k)^{1/2}\Big[\int_{U_0} g_0(x,u)^2
\mu_0(du)\Big]^{1/2},
 \eeqlb
where
 \beqlb\label{4.9}
\int_{V_k} |g_0(y_n,u)-g_0(x,u)| \mu_0(du)
 \le
\Big[\mu_0(V_k)\int_{U_0} |g_0(y_n,u)-g_0(x,u)|^2 \mu_0(du)\Big]^{1/2}.
 \eeqlb
By letting $n\to \infty$ and $k\to \infty$ in (\ref{4.8}) and using
condition (3.b) one can see that
 \beqlb\label{4.10}
\lim_{n\to \infty} \int_{V_n} D_{g_0(y_n,u)}f(x_n)\mu_0(du)
 =
\int_{U_0} D_{g_0(x,u)}f(x) \mu_0(du).
 \eeqlb
Similarly, for $n\ge k$ we have
 \beqnn
\ar\ar\Big|\int_{W_n} \itDelta_{g_1(x_n,u)}f(x_n)\mu_0(du) - \int_{U_2}
\itDelta_{g_1(x,u)}f(x) \mu_1(du)\Big| \cr
 \ar\ar\qquad
\le \|f^\prime\|\int_{U_2} |g_1(x_n,u)-g_1(x,u)| \mu_1(du) +
2\|f^\prime\|\varepsilon_k \cr
 \ar\ar\qquad\quad
+ \Big[\|f^\prime\||x_n-x| + |f(x_n)-f(x)|\Big]\mu_1(W_k).
 \eeqnn
Then letting $n\to \infty$ and $k\to \infty$ and using condition (3.a) one
sees
 \beqlb\label{4.11}
\lim_{n\to \infty}\int_{W_n} \itDelta_{g_1(x_n,u)}f(x_n)\mu_0(du)
 =
\int_{U_2} \itDelta_{g_1(x,u)}f(x) \mu_1(du).
 \eeqlb
In view of (\ref{4.10}) and (\ref{4.11}), it is obvious that $A_nf(x_n)\to
Af(x)$ as $n\to \infty$. \qed

\bproposition\label{t4.3} Suppose that conditions (4.a) and (3.a,b) are
satisfied. Then there exists a weak solution to (\ref{2.2}). \eproposition

\proof Following the proof of Lemma~4.3 in \cite{FuL10} it is easy to show
that $\{x_n(t): t\ge0\}$ is a tight sequence in the Skorokhod space
$D([0,\infty), \mbb{R})$. Then there is a subsequence $\{x_{n_k}(t):
t\ge0\}$ that converges to some process $\{x(t): t\ge0\}$ in distribution
on $D([0,\infty), \mbb{R})$. By the Skorokhod representation theorem, we
may assume those processes are defined on the same probability space and
$\{x_{n_k}(t): t\ge0\}$ converges to $\{x(t): t\ge0\}$ almost surely in
$D([0,\infty), \mbb{R})$. Let $D(x) := \{t>0: \mbf{P}\{x(t-)=x(t)\}= 1\}$.
Then the set $[0,\infty) \setminus D(x)$ is at most countable; see, e.g.,
\cite[p.131]{EtK86}. It follows that $\lim_{k\to \infty} x_{n_k}(t) =
x(t)$ almost surely for every $t\in D(x)$; see, e.g., \cite[p.118]{EtK86}.
From (\ref{4.7}) and Lemma~\ref{t4.2} it follows that (\ref{4.2}) is a
locally bounded martingale. Then we get the result by
Proposition~\ref{t4.1}. \qed

\bproposition\label{t4.4} Suppose that conditions (4.a) and (3.a,c) hold
with: (i)~$c=1, \alpha_m=2, p_m=1/2$; or (ii)~$c<1, \alpha_m<2,
1-1/\alpha_m< p_m\le 1/2$. Then there exists a weak solution to
(\ref{2.2}). \eproposition

\proof In condition (3.c), we can obviously assume $f_m\le f_{m+1}$ for
all $m\ge 1$. Let $V_n = \{u\in U_0: f_n(u)\ge 1/n\}$. Then the conclusion
of Lemma~\ref{t4.2} remains true. The only necessary modification of the
proof is that now we consider $n\ge k\ge M$. Then $|x|, |x_n|\le M$
implies $|x|, |y_n|\le k$, so we can replace (\ref{4.9}) by
 \beqnn
\int_{V_k} |g_0(y_n,u)-g_0(x,u)| \mu_0(du)
 \ar\le\ar
|y_n-x|^{p_k}\int_{V_k} f_k(u) \mu_0(du) \cr
 \ar\le\ar
k|y_n-x|^{p_k}\int_{U_0} [f_k(u)\land f_k(u)^2] \mu_0(du).
 \eeqnn
Then the result follows as in the proof of Proposition~\ref{t4.3}. \qed


\section{Strong solutions}

\setcounter{equation}{0}

In this section, we prove the existence of the strong solution to
(\ref{2.1}). Let $(\sigma,b,g_0,g_1)$ be given as in the second section.
We assume the following linear growth condition on the coefficients:
 \bitemize

\itm[{\rm(5.a)}] there is a constant $K\ge 0$ such that
 \beqnn
\ar\ar\sigma(x)^2 + \int_{U_0} g_0(x,u)^2\mu_0(du) + \int_{U_2}
g_1(x,u)^2\mu_1(du) \cr
 \ar\ar\qquad
+\, b(x)^2 + \bigg(\int_{U_2} |g_1(x,u)| \mu_1(du)\bigg)^2\leq K(1+x^2),
\quad x\in \mbb{R}.
 \eeqnn

 \eitemize

\btheorem\label{t5.1} Suppose that conditions (5.a) and (3.a,b) are
satisfied. Then there is a pathwise unique strong solution to (\ref{2.1}).
\etheorem

\proof By Proposition~\ref{t4.3} for each integer $m\ge 1$ there is a weak
solution to
 \beqlb\label{5.1}
x(t)
 \ar=\ar
x(0) + \int_0^t \sigma(\chi_m(x(s)))dB(s) + \int_0^t b(\chi_m(x(s)))ds \nnm \\
 \ar\ar\qquad
+ \int_0^t\int_{U_0} g_0(\chi_m(x(s-)),u) \tilde{N}_0(ds,du) \nnm \\
 \ar\ar\qquad
+ \int_0^t\int_{U_2} \chi_m\circ g_1(\chi_m(x(s-)),u) N_1(ds,du).
 \eeqlb
The pathwise uniqueness for the equation follows from Theorem~\ref{t3.2}.
Then there is a unique strong solution $\{x_m(t): t\ge 0\}$ to
(\ref{5.1}); see, e.g., \cite[p.104]{Sit05}. Let $\tau_m = \inf\{t\ge0:
|x_m(t)|\ge m\}$. As in the proof of Proposition~3.4 in \cite{LiM11} it is
easy to get
 \beqnn
\mbf{E}\Big[1+\sup_{0\le s\le t}x_m(s\land\tau_m)^2\Big]
 \ar\le\ar
(1+6\mbf{E}[x(0)^2])\exp\{6K(4+t)t\}.
 \eeqnn
Then $\tau_m\to \infty$ as $m\to \infty$. Following the proof of
Proposition~2.2 in~\cite{FuL10} one can show there is a pathwise unique
strong solution to (\ref{2.2}). Then the result follows from
Proposition~\ref{t2.1}. \qed

\btheorem\label{t5.2} Let $\alpha_m$ be the number defined in (\ref{3.5}).
Suppose that conditions (5.a) and (3.a,c) hold with: (i)~$c=1, \alpha_m=2,
p_m=1/2$; or (ii)~$c<1, \alpha_m<2, 1-1/\alpha_m< p_m\le 1/2$. Then there
exists a pathwise unique strong solution to (\ref{2.1}). \etheorem

\proof Based on Proposition~\ref{t4.4}, this follows similarly as
Theorem~\ref{t5.1}. \qed


\section{Non-negative solutions}

\setcounter{equation}{0}

In this section, we derive some results on non-negative solutions of the
stochastic equation (\ref{2.1}). Let $(\sigma,b,g_0,g_1)$ be given as in
the second section. In addition, we assume:
 \bitemize

\itm $b(x)\ge 0$ and $\sigma(x) = 0$ for $x\le 0$;

\itm for every $u\in U_0$ we have $x+g_0(x,u)\ge 0$ if $x>0$ and
$g_0(x,u)=0$ if $x\le 0$;

\itm $x+g_1(x,u)\ge 0$ for $u\in U_1$ and $x\in \mbb{R}$.

 \eitemize
Then, by Proposition~2.1 in \cite{FuL10}, any solution of (\ref{5.1}) is
non-negative. By considering non-negative solutions, we can weaken the
linear growth condition of the parameters into the following:
 \bitemize

\itm[{\rm(6.a)}] there is a constant $K\ge 0$ such that
 \beqnn
b(x) + \int_{U_2} |g_1(x,u)| \mu_1(du) \le K(1+x), \qquad x\ge 0;
 \eeqnn

\itm[{\rm(6.b)}] there is a non-decreasing function $x\mapsto L(x)$
     on $\mbb{R}_+$ so that
 \beqnn
\sigma(x)^2 + \int_{U_0} g_0(x,u)^2 \mu_0(du)\le L(x), \qquad x\ge 0.
 \eeqnn

 \eitemize

\btheorem\label{t6.1} Suppose that conditions (6.a) and (3.a,b) are
satisfied. Then for any $x(0)\in \mbb{R}_+$ there is a pathwise unique
non-negative strong solution to (\ref{2.1}). \etheorem

\proof By conditions (6.a) and (3.b) one can show that the parameters of
(\ref{5.1}) satisfy condition (4.a). Then for each integer $m\ge 1$ there
is a non-negative weak solution to (\ref{5.1}) by Proposition~\ref{t4.3}.
The pathwise uniqueness for (\ref{5.1}) holds by Theorem~\ref{t3.2}, so
there is a unique non-negative strong solution to (\ref{5.1}). Then the
result follows as in the proof of Proposition~2.2 in \cite{FuL10}. \qed

\bcorollary\label{t6.2} {\rm (Dawson and Li \cite{DaL12})} Given $0\le
x(0)\le 1$ there is a pathwise unique strong solution $\{x(t): t\ge 0\}$
to (\ref{1.3}) such that $0\le x(t)\le 1$ for all $t\ge 0$. \ecorollary

\proof Observe that $q(x,r) = 0$ for $x\le 0$ and $x\ge 1$. For any $0\le
x,z,r\le 1$ we have
 \beqnn
0\le x+zq(x,r) = z1_{\{r\le x\}} + (1-z)x\le 1.
 \eeqnn
Then $0\le x(0)\le 1$ implies $0\le x(t)\le 1$ for all $t\ge 0$. The
function $x\mapsto x + q(x,r)$ is clearly non-decreasing and for any $0\le
x, y\le 1$,
 \beqnn
\int_0^1\nu(dz)\int_0^1 z^2|q(x,r) - q(y,r)|^2dr
 \ar=\ar
[|x-y| - (x-y)^2]\int_0^1z^2\nu(dz) \cr
 \ar\le\ar
|x-y|\int_0^1 z^2\nu(dz).
 \eeqnn
Then the result follows by Theorem~\ref{t6.1}. \qed

\bcorollary\label{t6.3} {\rm (D\"{o}ring and Barczy \cite{DoB11})} Given
$x(0)\ge 0$ there is a unique non-negative strong solution to (\ref{1.4}).
\ecorollary

\proof It is easy to see that $x\mapsto x + g(x,u,r)$ is a non-decreasing
function. For any $x, y\ge 0$ we have
 \beqnn
\ar\ar\int_0^\infty dr\int_0^\infty(g(x,u,r) - g(y,u,r))^2 \mu_0(du)\cr
 \ar\ar\qquad
= \int_0^\infty(1-e^{-u})^2\mu_0(du)\big[x+y - 2(x^{-1}\land
y^{-1})xy\big] \cr
 \ar\ar\qquad
= \int_0^\infty(1-e^{-u})^2\mu_0(du)|x-y|.
 \eeqnn
By Theorem~\ref{t6.1} there is a unique non-negative strong solution to
the equation. \qed

\btheorem\label{t6.4} Suppose that conditions (6.a,b) and (3.a,c) hold
with: (i)~$c=1, \alpha_m=2, p_m=1/2$; or (ii)~$c<1, \alpha_m<2,
1-1/\alpha_m< p_m\le 1/2$. Then there exists a pathwise unique
non-negative strong solution to (\ref{2.1}). \etheorem

\proof This follows similarly as Theorem~\ref{t6.1}. Here condition (6.b)
is used to guarantee condition (4.a) is satisfied by the parameters of
(\ref{5.1}). \qed

\bigskip

\textbf{Acknowledgements.} The authors want to thank Professor M.\ Barczy
for his careful reading of the paper and pointing out a number of typos.
We would also like to acknowledge the Laboratory of Mathematics and
Complex Systems (Ministry of Education, China) for providing us the
research facilities.


\end{document}